\documentclass[11pt,fleqn,twoside]{article}
\usepackage{amsfonts,latexsym}
\makeatletter
\newcommand{\prava}{\footnotesize\it
\begin{flushright}
\begin{minipage}{18cm}
Copyright \copyright 1998 by A. Sat{\i}r
\end{minipage}
\end{flushright}}

\newcommand{\name}[1]{\begin{flushleft}
                       \LARGE \bf #1
                       \end{flushleft}\vspace{-3mm}}

\newcommand{\Author}[1]{\begin{flushleft}
                       \it #1 \end{flushleft}}

\newcommand{\Adress}[1]{\begin{flushleft}
                       \it #1 \end{flushleft}}

\newcommand{\Date}[1]{\begin{flushleft}
                      \small  \it #1 \end{flushleft}}

\newcommand{\ehkol}{Author \ name}
\newcommand{\ohkol}{Article \ name}
\renewcommand{\@evenhead}{
\hspace*{-3pt}\raisebox{-15pt}[\headheight][0pt]{\vbox{\hbox to \textwidth
{\thepage \hfil \ehkol}\vskip4pt \hrule}}}
\renewcommand{\@oddhead}{
\hspace*{-3pt}\raisebox{-15pt}[\headheight][0pt]{\vbox{\hbox to \textwidth
{\ohkol \hfil \thepage}\vskip4pt\hrule}}}
\renewcommand{\@evenfoot}{}
\renewcommand{\@oddfoot}{}

\newcommand{\be}{\begin{equation}}
\newcommand{\ee}{\end{equation}}
\newcommand{\ba}{\hspace*{-5pt}\begin{array}}
\newcommand{\ea}{\end{array}}

\newcommand{\ds}{\displaystyle}
\makeatother

\begin{document}
\setcounter{page}{364}
\thispagestyle{empty}

\renewcommand{\ehkol}{A. Sat{\i}r}
\renewcommand{\ohkol}{Dif\/ferential Constraints Compatible with
Linearized Equations}

\begin{flushleft}
\footnotesize \sf
Journal of Nonlinear Mathematical Physics \qquad 1998, V.5, N~4,
\pageref{satir-fp}--\pageref{satir-lp}. \hfill {\sc Letter}
\end{flushleft}

\vspace{-5mm}

\renewcommand{\footnoterule}{}
{\renewcommand{\thefootnote}{}
 \footnote{\prava}

\name{Dif\/ferential Constraints Compatible \\
with Linearized Equations} \label{satir-fp}

\Author{Ahmet SATIR}

\Adress{Middle East Technical University, Department of Physics,\\
06531 Ankara, Turkey\\[1mm]
Correspondence address: \\
Koza Sokak 136/1, G.O.P. 06670 Ankara, Turkey \\
Fax: +90 312 437 8634}

\Date{Received January 15, 1998; Revised July 30, 1998;
Accepted July 31, 1998}

\begin{abstract}
\noindent
Dif\/ferential constraints compatible with the linearized equations
of partial dif\/ferential equations are examined. Recursion operators
are obtained by integrating the dif\/fe\-ren\-tial constraints.
\end{abstract}

\noindent
One of the standard ways for determining particular solutions to partial
dif\/ferential equations is to reduce them to ordinary dif\/ferential equations
which are easier to solve. The classical work of Lie about group-invariant
solutions generalizes well-known methods for f\/inding similarity solutions
and other basic reduction methods \cite{satir:olver}. Bluman and
Cole~\cite{satir:blu} proposed a generalization of Lie's method for
f\/inding group-invariant solutions, which they named the ``nonclassical''
method. In this approach, one replaces the condition for the
invariance of the given system of
dif\/ferential equations by the weaker condition for the invariance of the
combined system consisting of the original dif\/ferential equations along
with the equations requiring the group invariance of the solutions.
P.J.~Olver and P.~Rosenau proposed a generalization of the
nonclassical method \cite{satir:olver1,satir:olver2}.
They showed that many known reduction methods,
including the classical and nonclassical methods, partial invariance, and
separation of variables can be placed into a general framework. In their
formulation, the original system of partial dif\/ferential equations can be
enlarged by appending additional dif\/ferential constraints (side conditions),
such that the resulting overdetermined system of partial dif\/ferential
equations satisfy compatibility conditions.

This work discusses
dif\/ferential constraints compatible with the linearized
equations of partial dif\/ferential equations instead of the partial
dif\/ferential equations themselves. The relation between dif\/ferential
constraints and recursion operators are examined. For the type of equations
in the form $q_t=P(x,t,q,q_x,q_{xx})$ and $q_t=P(q,q_x,q_{xx},q_{xxx})$,
recursion operators are obtained by integrating the compatible dif\/ferential
constraints. A new type of integrable equations, which are generalizations
of the integrable equations of Fokas and Svinolupov, are given. Results are
also compared with Fokas' generalized symmetry~\cite{satir:fokasj} and
Mikhailov-Shabat-Sokolov's formal symmetry
approaches~\cite{satir:svi,satir:mss}.

We can describe the dif\/ferential constraint method for evolutionary
equations~\cite{satir:ben}
\begin{equation}
q_t=P(x,t,q,q_x,q_{xx},\ldots)  \label{satir:qt}
\end{equation}
in the following way. First we linearize the given dif\/ferential
equation. In other words, we replace $q$ (and its derivatives) in
(\ref{satir:qt}) by $q+\epsilon \Psi$
and dif\/ferentiate both sides of the resulting expression
with respect to $\epsilon $ and take the limit $\epsilon \rightarrow
0$, i.e.,
\begin{equation}
\Psi _t=D_P(\Psi )  \label{satir:sym}
\end{equation}
where $D_P$ is the Fr{\'{e}}chet derivative \cite{satir:olver}.
The equation above
can also be written as
\begin{equation}
\Psi _t=\sum_{i=0}^N
P_i\Psi _i=\sum_{i=0}^N\frac{\partial P}{\partial q_i}\Psi _i
\end{equation}
where $N$ is the order of dif\/ferential equation, $q_0=q$,
$q_1=q_x$, $q_2=q_{xx},$
$\Psi _0=\Psi$, $\Psi _1=\Psi _x$, $\Psi _2=\Psi _{xx}$, and so on. In the
classical symmetry approach, (\ref{satir:sym}) is the main equation,
where
$\Psi $ is
the symmetry of the dif\/ferential equation, which is a function of
$x$, $t$, $q_i$.

The compatible dif\/ferential constraint is
\begin{equation}
H\Psi =0,  \label{satir:h}
\end{equation}
where $H$ depends on $x$, $t$, $q_i$. If its order
(highest derivative in $H$) is $N$, then (\ref{satir:h}) may be written as
\begin{equation}
\Psi _N=\sum_{i=0}^{(N-1)}A_i\Psi _i,  \label{satir:eigen}
\end{equation}
where $A_0,A_1,\ldots ,A_{N-1}$ are functions of $x$, $t$, $q_i$.
Compatibility of (\ref{satir:eigen}) and (\ref{satir:sym}) is given by
\begin{equation}
\Psi _{N,t}-\Psi_{t,N}=0.
\end{equation}
Using (\ref{satir:sym}) and (\ref{satir:eigen}), rewlation (6) leads to
\begin{equation}
\sum_{i=0}^{N-1}\Psi _{i}W_i=0.
\end{equation}
Letting
\begin{equation}
W_i=0,  \label{satir:u}
\end{equation}
we obtain a system of partial dif\/ferential equations among $P_i$,
$A_i$ and
their partial derivatives. The solution of this system will determine the
dif\/ferential constraint (\ref{satir:h}), which can be integrated to give
\begin{equation}
\Phi \Psi =0,  \label{satir:rec11}
\end{equation}
where $\Phi $ is the recursion operator.

Let us now consider dif\/ferential equations of the following form
\begin{equation}
q_t=f(x,t,q,q_x,q_{xx}).  \label{satir:q1}
\end{equation}
The linearized form of (10) can be given as
\begin{equation}
\Psi_t=\gamma \Psi _{xx}+\alpha \Psi _x+\beta \Psi,  \label{satir:q2}
\end{equation}
where $\alpha $, $\beta $, $\gamma $ are functions of $x$, $t$, $q$, $q_x$,
$q_{xx} $. We assume a dif\/ferential constraint having the same order as
(\ref{satir:q1}) in the form
\begin{equation}
\Psi _{xx}=A\Psi _x+B\Psi,   \label{satir:ei1}
\end{equation}
where $A$ and $B$ are functions of $x$, $t$, $q$, $q_x$, $q_{xx}$. The
compatibility condition between (\ref{satir:q2}) and
(\ref{satir:ei1}) will give the following evolution equations
\be
\ba{l}
A_t =\alpha _{xx}+\alpha _xA+\gamma _{xx}A+2\gamma _xA_x+\gamma
_xA^2+2\gamma _xB  \label{satir:a}
\vspace{2mm}\\
\qquad +2\beta _x+A_{xx}\gamma +A_x\alpha +2A_x\gamma A+2B_x\gamma  ,
\ea
\ee
\be
\ba{l}
B_t =2\alpha _xB+\gamma _{xx}B+2\gamma _xB_x+\gamma _xAB+\beta _{xx}
\label{satir:b}
\vspace{2mm}\\
\qquad -\beta _xA+2A_x\gamma B+B_{xx}\gamma +B_x\alpha .
\ea
\ee
The solutions of the system (\ref{satir:a})--(\ref{satir:b}) can be
given with the linearized equation
\begin{equation}
\Psi _t=\eta \Psi _{xx}+\left[\frac{2\eta r_{qq}}{r_q}q_x+2\eta r+\eta _1
\right]\Psi_x+
\left[\frac{\eta \left(r_{qqq}r_q-r_{qq}^2\right)}{r_q^2}q_x^2+2\eta
r_qq_x\right]\Psi,
\label{satir:lin}
\end{equation}
and compatible dif\/ferential constraint
\begin{equation}
\Psi_{xx}=\left[\frac{q_{xx}}{q_x}-\frac{q_xr_{qq}}{r_q}-r\right]\Psi_x
+\left[\frac{rq_{xx}}{q_x}-\frac{\left(r_{qqq}r_q-r_{qq}^2\right)q_x^2}
{r_q^2}-2r_qq_x-r_x\right]\Psi .
\label{satir:con2}
\end{equation}
Here $\eta $, $\eta _1$ are function of $x$, $t$, and $r$ is a
function of $x$, $t$, $q$, such that they satisfy
\begin{equation}
2r_x\eta +\eta _xr=0,\qquad 2r_t\eta +\eta _tr=0  \label{satir:cond1}
\end{equation}
\begin{equation}  \label{satir:cond2}
r_t-\frac r{4\eta }\left(-2\eta \eta _{xx}+\eta _x^2-2\eta _1\eta
_x+4\eta \eta _{1,x}\right)=0.
\end{equation}
It is interesting to note that in this classif\/ication the explicit time
dependence is crucial for obtaining equations
(\ref{satir:cond1})--(\ref{satir:cond2}),
which do not have any derivative with respect to $q$, although $r$ depend
explicitly on $q $. The dif\/ferential constraint (\ref{satir:con2}) can be
integrated to give
\begin{equation}
\Phi =D+\frac{r_{qq}q_x}{r_q}+r+q_xD^{-1}r_q+q_xD^{-1}\left(\frac{
r_{qqx}r_q-r_{qq}r_{qx}}{r_q^2}\right),  \label{satir:f2rec}
\end{equation}
where $\ds \left(D^{-1}f\right)(x)=\int_{-\infty }^xf(\zeta
)\,d\zeta$. Because of the last term in (\ref{satir:f2rec}) and
the explicit dependence of $x$ and $t$ on $r$, $\Phi $, as given by 
(\ref{satir:f2rec}), is new recursion operator which is a
generalization of the recursion operators obtained by Fokas \cite{satir:fokasj}
and Ibragimov-Shabat \cite{satir:ib}. The new integrable equation
takes the following form
\begin{equation}
q_t=\eta q_{xx}+\frac{\eta r_{qq}q_x^2}{r_q}+(2\eta r+\eta _1)q_x.
\label{satir:eq2}
\end{equation}
Although (20) looks like the integrable equations classif\/ied
by Fokas~\cite{satir:fokasj}, they are dif\/ferent, because of the 
explicit $x$, $t$ dependence of $r$, $\eta $ and $\eta_{1} $, which
are subject to conditions (\ref{satir:cond1})--(\ref{satir:cond2}).
Using transformation $q(x,t)\rightarrow
r(q(x,t),x,t)$~\cite{satir:fokasj}
one obtains
\begin{equation}
r_t=\eta r_{xx}+(2\eta r+\eta _1)r_x-\eta (r_{q,x}q_x+r_{,xx})-(2\eta r+\eta
_1)r_{,x}+r_{,t},  \label{satir:eq4}
\end{equation}
where $r_{,x}$, $r_{,t}$ denotes partial derivatives with respect to
$x$, $t$. However, in the limit $r(q,x,t)\rightarrow r(q)$, we
recover the equations classif\/ied by Fokas~\cite{satir:fokasj},
Ibragimov-Shabat~\cite{satir:ib} and Olver~\cite{satir:olv}.
Moreover, equations~(\ref{satir:eq4}) are more general than the once
obtained by
Svinolupov~\cite{satir:svi}, since he analyzed the case
$q_t=F(x,q,q_x,q_{xx})$, in which there
is no explicit time dependence. It is interesting to note that in the
limit $r(q,x,t)\rightarrow r(q,t)$ (using the conditions
(\ref{satir:cond1})--(\ref{satir:cond2})), equation~(\ref{satir:eq2})
will give Svinolupov's equation~\cite{satir:svi}
\begin{equation}  \label{satir:svi}
q_t=\eta ^{\prime }q_{xx}+\frac{\eta ^{\prime }r_{qq}q_x^2}{r_q}+\left(2\eta
^{\prime }r-\frac{\eta _t^{\prime }}{2\eta ^{\prime }}x+\epsilon \right)q_x,
\label{satir:eq3}
\end{equation}
where $\eta ^{\prime }$ and $\epsilon $ are functions of $t$, and $r$ is a
function of $q$ and $t$.

Equations (\ref{satir:cond1}) can be integrated to give $\ds \eta
=\frac a{r^2}$, where $a=a(q)$. Substitution of $\eta $
in~(\ref{satir:cond2}) will give
\begin{equation}
ar_{xx}r-2ar_x^2+r_x\eta _1r^3+\eta _{1,x}r^4-r_tr^3=0 .  \label{satir:devam}
\end{equation}
One particular solution can be given by solving the last two terms in the
above equation. This case ($r(q,t))$ corresponds to
equation~(\ref{satir:svi}). 
Equation~(\ref{satir:devam}) is a coupled second order partial dif\/ferential
equation and the general solution can be given by the method of
characteristics~\cite{satir:sneddon}. It is hoped that the general
solution will appear in a future work.

Next, we consider dif\/ferential equations of the following general
form,
which includes the KdV equation,
\begin{equation}
q_t=P(q,q_x,q_{xx},q_{xxx}).  \label{satir:q3}
\end{equation}
The linearization of (24) takes the form
\begin{equation}
\Psi _t=\alpha \Psi _{xxx}+\beta \Psi _{xx}+\gamma \Psi _x+\delta \Psi,
\label{satir:g1lin}
\end{equation}
where $\alpha $, $\beta $, $\gamma $, $\delta $ are functions of $q$,
$q_x$, $q_{xx}$, $q_{xxx}$. We consider the dif\/ferential constraint
having the same order as (\ref{satir:q3}), i.e.,
\begin{equation}
\Psi_{xxx}=A\Psi _{xx}+B\Psi _x+C\Psi   \label{satir:g1ei}
\end{equation}
where $A$, $B$ and $C$ are functions of $q$, $q_x$, $q_{xx}$, $q_{xxx}$. The
compatibility will determine algebraic equations
\begin{equation}
B=-\frac{2\alpha }{3\eta },\qquad C=-\frac 1{3\eta }(\alpha _x-2\alpha
A+2\beta )
\end{equation}
and evolution equations 
\be
A_t =\alpha _xA+\beta _x+A_{xxx}\eta +3A_{xx}\eta A
+3A_x^2\eta +A_x\alpha +3A_x\eta A^2,
\ee
\be
\ba{l}
\ds \alpha _t =\frac 12(2\alpha _{xxx}\eta -3\alpha _{xx}\eta A-3\alpha
_xA_x\eta +2\alpha _x\alpha  \vspace{2mm}\\
\ds \qquad -3\beta _{xx}\eta +6\beta _x\eta A+6A_x\beta \eta
+3\varepsilon _t\eta )
\ea
\ee
\be
\ba{l}
\ds \beta _t =\frac 14(3\alpha _{xxx}\eta A-6\alpha _{xx}\eta A^2-3\alpha
_xA_{xx}\eta  -12\alpha _xA_x\eta A+4\alpha _x\beta
\vspace{3mm}\\
\ds \qquad +\beta _{xxx}\eta -6\beta _{xx}\eta
A+4\beta _x\alpha +12\beta _x\eta A^2+6A_{xx}\beta \eta +24A_x\beta
\eta A),
\ea
\ee
where $\eta $ is constant and $\epsilon $ is a function of $t$.

The linearized form of the f\/irst class is
\begin{equation}
\Psi_{t}=\eta \Psi{xxx} + \left[\frac{\rho_{1}}{2}q_{x}^{2} + \rho
\right] \Psi_{x} + \rho_{q} q_{x} \Psi,
\end{equation}
with compatible dif\/ferential constraint
\begin{equation}  \label{satir:rec1}
\Psi_{xxx}=\frac{q_{xx}}{q_{x}} \Psi_{xx} - \left[\frac{\rho_{1}}{3 \eta}
q_{x}^{2} + \frac{2 \rho}{3 \eta} \right] \Psi_{x} -\frac{ -2 \rho q_{xx} + 3
\rho_{q} q_{x}^{2}} {3 \eta q_{x} } \Psi,
\end{equation}
where $\rho_{1} $ is constant and $\rho $ is a function of $q$ with the
condition
\begin{equation}
\rho_{qqq}+\frac{4 \rho_{1}}{3 \eta} \rho_{q}=0.
\end{equation}
The recursion operator can be obtained by integrating
(\ref{satir:rec1}).
This leads to
\begin{equation}  \label{satir:kdvr1}
\Phi = D^{2} + \frac{2 \rho }{ 3 \eta } + \frac{\rho_{1}}{3 \eta} q_{x}^{2}
- \frac{\rho_{1}}{3 \eta} q_{x} D^{-1} q_{xx} + \frac{q_{x}}{ \eta} D^{-1}
\rho_{q}.
\end{equation}
The integrable equation is in the form
\begin{equation}
q_t=\eta q_{xxx}+\frac{\rho _1}6q_x^3+\rho q_x.  \label{satir:kdv1}
\end{equation}

The second class is given by the linearized equation
\begin{equation}
\Psi _t=\eta \Psi_{xxx}+\left[\epsilon _1q_x^2+2\epsilon _1\epsilon
_2q_x+2\epsilon _3\right]\Psi
\end{equation}
with compatible dif\/ferential constraint
\be \label{satir:rec2}
\Psi_{xxx}=\frac{q_{xx}}{q_x+\epsilon _2}\Psi _{xx} -\left[\frac{\epsilon
_1}{3\eta }q_x^2+\frac{2\epsilon _1\epsilon _2}{3\eta }q_x +
\frac{2\epsilon _3}{ 3\eta }\right]\Psi _x
-\frac{\left(\epsilon _2^2\epsilon _1-2\epsilon
_3\right)q_{xx}}{3\eta (q_x+\epsilon _2)} \Psi,
\ee
where $\epsilon _1$, $\epsilon _2$, $\epsilon _3$ are constants. The recursion
operator can be obtained by integrating (\ref{satir:rec2}), which gives
\begin{equation}  \label{satir:kdvr2}
\Phi =D^2+\frac{2\epsilon _3}{3\eta }+\frac{\epsilon _2}{3\eta }q_x^2+\frac{%
2\epsilon _1\epsilon _2}{3\eta }q_x-\frac{\epsilon _1}{3\eta }
(q_x+\epsilon_2) D^{-1}q_{xx}.
\end{equation}
The integrable equation is in the form
\begin{equation}  \label{satir:kdv2}
q_t=\eta q_{xxx}+\frac{\epsilon _1}6q_x^3+\frac{\epsilon _1\epsilon _2}%
2q_x^2+\epsilon _3q_x.
\end{equation}
Equations (\ref{satir:kdv1}), (\ref{satir:kdv2}) were classif\/ied by
Fokas \cite{satir:fokasj},
Ibragimov-Shabat \cite{satir:ib}. We obtained the recursion ope\-ra\-tors
(\ref{satir:kdvr1}), (\ref{satir:kdvr2})
by integrating the dif\/ferential constraints.

The third type of equation is given by the linearized equation
\begin{equation}
\Psi _t=\lambda _5\Psi _{xxx}+\frac 32\lambda _3\lambda _5\Psi _{xx}+\frac
32\lambda _4\lambda _5\Psi _x+\lambda _6\Psi,
\end{equation}
and dif\/ferential constraint
\be
\ba{l}
\ds \Psi _{xxx} =\left[\frac{q_{xxx}+\lambda _1q_{xx}+\lambda
_2q_x}{q_{xx}+ \lambda _1q_x+\lambda _2q}-\lambda _3\right]\Psi _{xx}
\vspace{3mm}\\
\ds \qquad+\left[\frac{\lambda _{3(}q_{xxx}+\lambda _1q_{xx}+\lambda
_2q_x)}{q_{xx}+\lambda _1q_x+\lambda _2q}-\lambda _4\right]\Psi _x
+\frac{\lambda _{4(}q_{xxx}+\lambda _1q_{xx}+\lambda _2q_x)}{%
q_{xx}+\lambda _1q_x+\lambda _2q}\Psi .
\ea
\ee
The recursion operator is
\begin{equation}
\Phi =D^2+\lambda _3D+\lambda _4 .
\end{equation}
The integrable equation is given as
\begin{equation}
q_t=\lambda _5q_{xxx}+\frac 32\lambda _3\lambda _5q_{xx}+\frac 32\lambda
_4\lambda _5q_x+\lambda _6q,
\end{equation}
where $\lambda _1$, $\lambda _2$, $\lambda _3$, $\lambda _4$,
$\lambda _5$ and $\lambda _6$ are constants. Note that Rabelo and Tanenblat
also obtained linear equation using the classif\/ication method of
pseudo-spherical surfaces with Gaussian curvature $(-1)$~\cite{satir:keti}.

According to Fokas, an integrable equation has inf\/initely many generalized
symmetries which are the solutions of (\ref{satir:sym}). The existence of
generalized symmetry manifest itself by the existence of an
admissible Lie-B{\"{a}}cklund operator. The existence of inf\/initely
many symmetries is expressed by the existence of a recursion
operator. There is a close relationship between Lie-B{\"{a}}cklund
operator and linearized equation \cite{satir:fokasj}. Because, if
Fokas' admissible Lie-B{\"{a}}cklund operator is applied on the
evolution equation, we obtain the Fr{\'{e}}chet derivative of our
dif\/ferential equation under consideration or linearized form of our
dif\/ferential equation. Recursion operators, in our method, are obtained
by the integration of the compatible dif\/ferential constraints.

Let us brief\/ly recall Olver-Fokas symmetry approach and the
dif\/ferential constraint test:
\begin{enumerate}
\item[] {\it Olver-Fokas symmetry test:} The equation $q_t=P[q]$
is integrable if there exists inf\/initely many non-Lie point symmetries, or
equivalently, one non-Lie point symmetry and a recursion operator.
The recursion
operator and the time-independent part of the linearized dif\/ferential
equation form
a Lax pair $\Re _t+[\Re ,D_P]=0$.

\item[] {\it Dif\/ferential constraint test:}
The equation $q_t=P[q]$ is integrable if there exists a dif\/ferential
constraint $H\Psi =0$, such that it is compatible with the linearized
equation $\Psi _t=D_P(\Psi )$. The compatibility condition is
$H_t+[H,D_P]=0$.
\end{enumerate}

There are several methods to examine the integrability of nonlinear
partial dif\/ferential equations \cite{satir:fokass}, although in two
dimensions most of these methods imply each other. Inf\/inite sets of
conservation laws, inf\/inite number of symmetries, and the bi-Hamiltonian
structure, are a number of remarkable properties, to name but three.
The recursion operator plays an important role in the formulation of
these recursive properties. Firstly, the fa\-mi\-ly of integrable
equations can be written in a compact form using recursion operators.
The other role of recursion operator is associated with the
Hamiltonian treatment of integrable equations. Recursion operators
determine Hamiltonian structures through certain Poisson
brackets~\cite{satir:fordy}. Further analysis of the recursion
operators leads to the concept of bi-Hamiltonian structures, which can
be given by factorizing the recursion operator.

The main idea in this work is to give a new def\/inition of
integrability. A partial dif\/ferential equation is integrable if its
linearized equation is compatible with a dif\/ferential constraint.
Using dif\/ferential constraints compatible with the linearized
equations, we also give an answer to the question ``The deep
connection between the direct reduction and recursion operators'' of
Olver~\cite{satir:olv}.

\subsection*{Acknowledgments}

I would like to thank Norbert Euler and the referee for
their constructive comments.

\label{satir-lp}


\begin{thebibliography}{99}
\footnotesize

\bibitem{satir:olver}  Olver P.J.,  Applications of Lie Groups to
Dif\/ferential Equations,  Springer-Verlag, Berlin, 1993.

\bibitem{satir:blu}  Bluman G.W. and Cole J.D., The General Similarity
Solutions
of the Heat Equation, {\it  J. Math. Mech.}, 1969, V.18, 1025--1042.


\bibitem{satir:olver1}  Olver P.J. and  Rosenau P.,
Group Invariant Solutions of
Dif\/ferential Equations, {\it SAIM J. Appl. Math.}, 1987, V.47,
263--278.

\bibitem{satir:olver2}   Olver P.J. and  Rosenau P.,
The Construction of Special
Solutions to Partial Dif\/ferential Equations, {\it  Phys. Lett.  A}, 1986,
V.11, 107--112.



\bibitem{satir:fokasj}   Fokas A.J.,
A Symmetry Approach to Exactly Solvable
Evolution Equations, {\it  J. Math. Phys.}, 1980, V.21,  1318--1325.


\bibitem{satir:svi}   Svinolupov S.I., 1985 Second-Order Evolution
Equations with
Symmetries, {\it  Russian Math. Surveys}, 1980, V.40,  241--242.


\bibitem{satir:mss}  Mikhailov A.V.,  Shabat A.B. and Sokolov V.V., in:
 What is Integrability, ed. V.E. Zakharov, Springer-Verlag, Berlin, 1991.



\bibitem{satir:olv} Olver P.J.,
Direct Reduction and Dif\/ferential
Constraints, {\it Proc.R.Soc.Lon. A}, 1994, V.444, 509--523.


\bibitem{satir:ben}  Sat{\i}r A., Dif\/ferential Constraints,
Recursion Operators and
Logical Integrability, {\it Inter. J. Theor. Phys.}, 1997, V.10, 2099--2105.


\bibitem{satir:ib}   Ibragimov N.Kh. and  Shabat A.B.,
Evolutionary Equations with
Nontrivial Lie-B\"{a}cklund Group, {\it  Funcl. Anal. Appl.},
1980, V.14, 19--28.


\bibitem{satir:sneddon}   Sneddon I.,
Elements of Partial Dif\/ferential Equations,
McGraw-Hill Kogakuska Ltd., 1957.


\bibitem{satir:keti}  Rabelo M.L. and  Tanenblat K.,
A Classif\/ication of
Pseudospherical Surface Equations of Type
$u_t=u_{xxx}+G(u,u_x,u_{xx})$, {\it  J. Math. Phys.}, 1992,
V.33, 537--549.


\bibitem{satir:fokass}  Fokas A.S., Symmetry and Integrability,
{\it Stud. Appl. Math.},  1987, V.77, 253--299.


\bibitem{satir:fordy}  Fordy A.P. and  Gibbons J.,
Factorization of Operators I. Miura Transformations,
{\it J. Math. Phys.}, 1980, V.21, 2508--2510.

\end{thebibliography}
\end{document}